\documentclass[11pt,a4paper]{amsart}
\usepackage{amsmath,amsthm,amssymb}
\usepackage{comment}
\newcommand{\mychoice}[3]{#1
}
\mychoice{
\newcommand{\plabel}[1]{ \label{#1}}
\newcommand{\gbibitem}[1]{ \bibitem{#1}}
\newcommand{\snewpage}{}

\specialcomment{commentx}{}{}\excludecomment{commentx}
\specialcomment{commenty}{}{}\excludecomment{commenty}
\specialcomment{commentz}{}{}

}{
\usepackage{xcolor}
\newcommand{\plabel}[1]{ \label{#1}\rlap{\smash{${}^{^{[#1]}}$}}}
\newcommand{\gbibitem}[1]{ \bibitem{#1}\rlap{\smash{${}^{^{[#1]}}$}}}
\newcommand{\snewpage}{\newpage}

\newenvironment{commentx}{\color{magenta} }{\color{black} }
\newenvironment{commenty}{\color{blue} }{\color{black} }

}{
\usepackage{xcolor}
\newcommand{\plabel}[1]{ \label{#1}}
\newcommand{\gbibitem}[1]{ \bibitem{#1}}
\newcommand{\snewpage}{}

\specialcomment{commentx}{}{}\excludecomment{commentx}
\newenvironment{commenty}{  }{}
\specialcomment{commentz}{}{}\excludecomment{commentz}

}

\DeclareMathOperator{\des}{des}
\DeclareMathOperator{\asc}{asc}

\DeclareMathOperator{\nonass}{n-a}

\DeclareMathOperator{\Lie}{Lie}

\DeclareMathOperator{\BCH}{BCH}
\DeclareMathOperator{\bch}{bch}
\DeclareMathOperator{\id}{id}

\newcommand{\bo}{\boldsymbol}
\newcommand{\lldots}{...}
\newcommand{\botimes}{{\textstyle\bigotimes}}
\newcommand{\bodot}{{\textstyle\bigodot}}
\newcommand{\leaveout}[1]{}

\newcommand{\ass}{\mathrm{assoc}}

\theoremstyle{definition}
\newtheorem{point}{}[section]
\newtheorem{remark}[point]{Remark}

\newtheorem{defi}[point]{Definition}

\theoremstyle{plain}
\newtheorem{prop}[point]{Proposition}

\newtheorem{propdef}[point]{Proposition/Definition}

\newtheorem{lemma}[point]{Lemma}
\newtheorem{theorem}[point]{Theorem}

\newtheorem{cor}[point]{Corollary}
\newcommand{\qedremark}{  \renewcommand{\qedsymbol}{$\triangle$} \qed \renewcommand{\qedsymbol}{$\Box$}}
\newcommand{\qedno}{\renewcommand{\qedsymbol}{}}
\newcommand{\marginextend}[1]{ \addtolength{\oddsidemargin}{-#1}  \addtolength{\evensidemargin}{-#1}\addtolength{\textwidth}{#1}\addtolength{\textwidth}{#1}}
\newcommand{\updownextend}[1]{ \addtolength{\topmargin}{-#1}  \addtolength{\textheight}{#1}
\addtolength{\textheight}{#1}}
\marginextend{1.25cm}
\updownextend{0cm}
\allowdisplaybreaks[3]
\begin{document}
\title[The Poincar\'e--Birkhoff--Witt theorem and Dynkin--Magnus commutators]{The symmetric Poincar\'e--Birkhoff--Witt theorem and Dynkin--Magnus commutators}
\date{\today}
\author{Gyula Lakos}
\email{lakos@renyi.hu}
\address{Alfréd Rényi Institute of Mathematics}
\begin{abstract}
The objective of this paper is to give alternative proofs for the symmetric Poincar\'e--Birkhoff--Witt theorem
 utilizing the Magnus recursion formulae or Dynkin's non-com\-mu\-ta\-ti\-ve polynomial comparison method
 and simple universal algebraic principles.
As an application of these principles, a theorem of  Nouaz\'e--Revoy type is also obtained.
\end{abstract}
\maketitle

%\snewpage

\textsl{To the 70th anniversary of}
Wilhelm Magnus:
On the exponential solution of differential equations for a linear operator.
\textit{Comm. Pure Appl. Math.}  \textbf{7}  (1954), pp. 649--673.

\textsl{To the 75th anniversary of}
E. B. Dynkin:
On the representation by means of commutators of the series $\log(e^xe^y)$ for noncommutative $x$ and $y$. (In Russian.)
\textit{Mat. Sb. (N.S.)}, \textbf{25(67)} (1949), pp. 155--162.
\section*{Introduction}\plabel{sec:introSymm}
The objective of this paper is to give alternative proofs for the symmetric Poincar\'e--Birkhoff--Witt theorem
 utilizing the Magnus recursion formulae  or Dynkin's non-com\-mu\-ta\-ti\-ve polynomial comparison method
 and simple universal algebraic principles.

\textbf{Historical context.}
The Baker--Campbell--Hausdorff expansion and the Poincar\'e--Birkhoff--Witt theorem
 have common roots as both of them originate from the study of Lie's Third Theorem.

The early history of the BCH expansion is narrated expertly in Achilles, Bonfiglioli \cite{AB};
 therefore we highlight just some main points:
Campbell \cite{Cam} demonstrated that the BCH expansion can formally be represented by a Lie series.
A simple formula for the generation of BCH expansion was given by Baker \cite{Bak}.
The resulting algebraic machinery was further clarified by Hausdorff \cite{Hau}, who also dealt with convergence.
At that point, the geometric-analytic content was still dominant in a certain sense:
Even if the symbolic nature of the Baker--Hausdorff formula is evident,
 the formal substitution operations used were motivated by vector fields.
The next wave of developments occurred after the corresponding algebraic machinery got clarified.
In that a crucial step was the Poincaré--Birkhoff--Witt theorem (see it later).
Simultaneously and consequently, an understanding of the nature of free Lie algebras relative to their
 (non-commutative polynomial) enveloping algebras was obtained, see Magnus \cite{MM} and Witt \cite{W}.
Thus the study of $\log((\exp X)(\exp Y))$ as an associative algebraically induced object became sufficient.
In that setting, Dynkin \cite{Dy} shortcut the earlier algorithmic and convergence arguments
 by the use of his polynomial commutator expansion comparison method
 which later became abstracted as the Dynkin--Specht--Wever lemma (cf.~Specht \cite{Sp}, Wever \cite{We}).
A need for a qualitative argument (for having Lie series), however, remained.
This was elegantly solved by the application of the criterion of Friedrichs,
 see Magnus \cite{M} and Finkelstein \cite{Fin} (cf.~Reutenauer \cite{R} for further references).
Meanwhile, by Dynkin \cite{Dyy} an alternative approach to the Lie series problem was offered.
Another possible way to view the BCH expansion is via its generalization, the Magnus expansion of Magnus \cite{M}.
Ultimately, Dynkin \cite{Dyy} and Magnus \cite{M} not only algebraize the BCH expansion but transcend it by finer methods.

From the preceding it must be clear that how the work of Magnus influenced the study of the broader area of the BCH expansion
 by his study of free Lie algebras \cite{MM} and the Friedrichs criterion \cite{M};
 and also by the direct generalization, the Magnus expansion \cite{M}.
The universal algebraic principles of Dynkin \cite{Dy}, \cite{Dyy} are
also clearly relevant although his latter work is less known.
What might not be immediately clear is whether there is a relevance in the other direction,
 toward the PBW theorem, thus, in particular, free Lie algebras.
The comparison of Dynkin \cite{Dyy}, Magnus \cite{M}, Solomon \cite{SS}, Mielnik, Pleba\'nski \cite{MP} shows that there should be.
The goal of this paper is articulate this connection.
Altogether, one expects that the standard universal Lie algebraic constructions ``support each other''.
Arguments of this type were considered before; see Cartier \cite{Ct} and Bonfiglioli, Fulci \cite{BF}, Ch.~6.
A difference compared to them is that  we use the finer Dynkin--Magnus commutators instead of the BCH  series.
%\snewpage

\textbf{The Poincaré--Birkhott--Witt theorem.}
Assume that $K$ is a unital commutative ring, and $\mathfrak g$ is a $K$-module with a compatible Lie-ring structure;
i.~e.~$\mathfrak g$ is a Lie $K$-algebra (also called: Lie ring over $K$).
The universal enveloping algebra $\mathcal U\mathfrak g$ is the free $K$-algebra
$\mathrm F_K[\mathfrak g]\simeq \bigotimes\mathfrak g\equiv \bigoplus_{n=0}^\infty\botimes^n\mathfrak g$ factorized by the ideal $J\mathfrak g$ generated by the
elements $X\otimes Y-Y\otimes X-\boldsymbol[X,Y\boldsymbol]$, the tensor products are taken over $K$.
Let $\boldsymbol m:\bigotimes\mathfrak g\rightarrow\mathcal U\mathfrak g$ denote this canonical homomorphism.

(i) If $\mathfrak g$ is a sum of cyclic $K$-modules, then we can choose a basis $\{g_\alpha\,:\,\alpha\in A\}$, and an ordering $\leq$ of $A$.
Then let $\bigotimes_\leq\mathfrak g$ be the submodule of $\bigotimes\mathfrak g$ spanned by $g_{\alpha_1}\otimes \ldots\otimes g_{\alpha_n}$ with $\alpha_1\leq\ldots\leq\alpha_n$.
Then the ``basic'' version of the PBW theorem states that
\begin{equation*}
\boldsymbol m_\leq:\textstyle{\bigotimes_\leq\mathfrak g}\rightarrow\mathcal U\mathfrak g
\end{equation*}
(a restriction of $\boldsymbol m$) is an isomorphism of $K$-modules.
This is due to Poincar\'e \cite{P}: $K$ is a field, $\mathbb Q\subset K$, cf.~Ton-That, Tran \cite{TT};
 and Birkhoff \cite{B}, Witt \cite{W}: $K$ is a field, but their methods work more generally.
Jacobson \cite{J} writes up the somewhat terse account of Birkhoff \cite{B};
 but the result is also not much different from the ideas of Witt \cite{W}.

(ii) A variant due to Cohn \cite{CC} is as follows:
If $\mathbb Q\subset K$, then one can consider the submodule $\bigotimes_\Sigma\mathfrak g$ of $\bigotimes\mathfrak g$.
This  submodule can be interpreted either as the submodule of elements invariant under permutations 
 in the order of tensor product or as the span of the elements
$a_1\otimes_\Sigma \ldots\otimes_\Sigma a_n=\frac1{n!}\sum_{\sigma\in\Sigma_n}g_{\sigma(1)}\otimes\ldots\otimes g_{\sigma(n)}$.
Then the ``symmetric'' version of the PBW theorem states that
\begin{equation*}
\boldsymbol m_\Sigma:\textstyle{\bigotimes_\Sigma\mathfrak g}\rightarrow\mathcal U\mathfrak g
\end{equation*}
(a restriction of $\boldsymbol m$) is an isomorphism of $K$-modules.

If $\mathbb Q\subset K$ is a field (which is quite the most important case), then it easy to see that the statements of
 the basic and symmetric versions are equivalent (cf. Bourbaki \cite{BX}).
In fact, they are connected by the local formulation of the PBW theorem:

In general, the enveloping algebra is naturally filtered by $\mathcal U^n\mathfrak g=\bo m\left(\botimes^{\leq n}\mathfrak g\right)$, and the
construction implies the existence of natural (surjective) maps
\[\bo m^{(n)}: \bodot^n\mathfrak g\rightarrow \mathcal U^n\mathfrak g/\mathcal U^{n-1}\mathfrak g\]
(factors of $\boldsymbol m$).
The local form of the Poincar\'e--Birkhoff--Witt theorem, whenever it holds, states  that the maps $\bo m^{(n)}$ are isomorphisms.
This holds is the cases (i) and (ii) above, i.~e.~when
$\mathfrak g$ is a sum of cyclic $K$-modules or if $\mathbb Q\subset K$, respectively.
The local statements are seemingly weaker is than the global ones but it is easy to see equivalence
(having $ \boldsymbol m_\leq$ and $\boldsymbol m_\Sigma$ defined). 
The validity of local PBW theorem, however, can be asked for any Lie
$K$-algebra $\mathfrak g$. Beyond cases (i) and (ii) it is known to hold if

(i') $K$ is a principal ideal domain (Lazard \cite{L}) or just a Dedekind domain (Cartier \cite{C});
 also see Higgins \cite{H} for further results in this direction;

(ii') $\frac12\in K$ but $[\mathfrak g,[\mathfrak g,\mathfrak g]]=0$ (Nouaz\'e, Revoy \cite{NR}).
\\But there are counterexamples
(see \v{S}ir\v{s}ov \cite{S}, Cartier \cite{C}, Cohn \cite{CC}).
%The most general approach is of Higgins \cite{H}, cf.~Revoy \cite{Re}.
For a broader view on the PBW theorem, we quote Higgins \cite{H}, Grivel \cite{G}, Shepler, Witherspoon \cite{SW}.

With the exception of the last section, we will be  concerned with the symmetric version of the PBW theorem,
independently from the general context above.
One can prove the surjectivity of $\boldsymbol m_\Sigma$
 using the ``symmetric rearrangement procedure'' in $\mathcal U\mathfrak g$,
 i.~e.~symmetrizing in the formally top nonarranged degree term
 at the cost of generating lower order terms, and repeating the process in formally lower orders.
The injectivity of $\boldsymbol m_\Sigma$ is, however, nontrivial.

%\snewpage

\textbf{Dynkin--Magnus commutators.}
Higher Lie algebraic commutators over $\mathbb Q\subset K$, can be described
(according to Solomon \cite{SS}) as objects  $\mu_n^{\Lie}$ ($n\geq1$)
such that

\begin{enumerate}
\item[($\mu1$)] $\mu^{\Lie}_1(X_1)=X_1$;
\item[($\mu2$)] $\mu^{\Lie}_n(X_1,\ldots,X_n)$ is a linear combination of Lie-monomials where every
variable has multiplicity $1$;
\item[($\mu3$)] the identities
\begin{equation}
\mu^{\Lie}_n(\ldots_1,X_{k-1},X_k,\ldots_2)-\mu^{\Lie}_n(\ldots_1,X_k,X_{k-1},\ldots_2)
=\mu^{\Lie}_{n-1}(\ldots_1,\bo[X_{k-1},X_k\bo],\ldots_2)
\plabel{eq:muid}
\end{equation}
hold for $n\geq2$, $1<k\leq n$.
\end{enumerate}
Thus these are Lie-polynomials $\mu_n^{\mathrm{\Lie}}(X_1,\ldots,X_n)\in\mathrm F^{\Lie}_K[X_1,\ldots,X_n ]$,
 i.~e.~elements of the free Lie $K$-algebra.
(Note on terminology:
For us commutator polynomials are certain elements of the non-commutative polynomial algebras $\mathrm F_K[X_1,\ldots,X_n ]$,
 while Lie-polynomials are elements of free Lie algebras $\mathrm F^{\Lie}_K[X_1,\ldots,X_n ]$, which
 evaluate to commutator polynomials naturally.
We also emphasize that $\mathrm F^{\Lie}_K[X_1,\ldots,X_n ]$ requires no particular ``construction'',
 it exists by universal algebraic reasons.)

The unicity of the $\mu^{\Lie}_n$ is relatively easy,
 but their existence, although they are sufficient to be constructed only over $\mathbb Q$, is less so.
Depending on viewpoint, these higher commutators can be called as Dynkin commutators (cf.~Dynkin \cite{Dyy}),
 Magnus commutators (cf.~Magnus \cite{M} and Mielnik, Pleba\'nski \cite{MP}),
 or first canonical projections (cf.~Solomon \cite{SS}).
In the view of Dynkin and Magnus, the higher commutators are the multilinear parts of the multivariable BCH expression.
In the view of Solomon, they can be considered as (representatives for) components of $(\boldsymbol m_\Sigma)^{-1}$.
As such, they seem to require some familiarity with either the BCH expansion or the PBW theorem itself.

Our point, however, is that
(a) $\mu_n^{\Lie}(X_1,\ldots,X_n)\in\mathrm F^{\Lie}_K[X_1,\ldots,X_n ]$ can be constructed directly,
(b) their existence is nearly synonymous to the symmetric PBW theorem,
(c) they can naturally be applied to approach the BCH / Magnus expansions.

\textbf{Outline of content.}
In Section \ref{sec:mu1}, we construct the $\mu_n^{\mathrm{\Lie}}$ using the idea of the Magnus recursion formulae.
In Section \ref{sec:mu2}, we construct the $\mu_n^{\mathrm{\Lie}}$ using ideas from Dynkin's noncommutative polynomial comparison method.
In Section \ref{sec:symm}, we prove the symmetric PBW theorem.
In Section \ref{sec:furthermu1}, further construction and comments related to $\mu_n^{\mathrm{\Lie}}$ are made.
In Section \ref{sec:udir}, some consequences of our viewpoint regarding the PBW theorem
 are discussed, a general theorem of  Nouaz\'e--Revoy type is obtained.

\textbf{Acknowledgements.}
 The author thanks Bal\'azs Csik\'os and M\'arton Nasz\'odi.
 The author is also grateful to David Yost for giving help regarding some of the literature.

\snewpage

\section{The existence of $\mu^{\Lie}$ I}\plabel{sec:mu1}
For practical reasons we will use left-iterated higher commutators
$\bo[X_1,\ldots, X_n\bo]_{\mathrm L}=\bo[ X_1,\bo[X_2,\ldots, \bo[X_{n-1},X_n\bo]\ldots\bo]\bo]$.
\begin{propdef}\plabel{pd:magnus}
There is a series of Lie-polynomials $ \mu^{\Lie}_n(X_1,\ldots,X_n)$, $n\geq1$, over $\mathbb Q$  such that
the following hold:

\begin{align}
&\mu^{\Lie}_n(X_1,\ldots,X_n)=\tag{L}\plabel{L}\\
&=\sum_{\substack{I_1\dot\cup\ldots \dot\cup I_s=\{2,\ldots,n\}\\I_k=\{i_{k,1},\ldots,i_{k,p_k}\}\neq
\emptyset \\i_{k,1}<\ldots<i_{k,p_k}}}\beta_s\cdot \bo[\mu^{\Lie}_{p_1}(X_{i_{1,1}},\ldots ,X_{i_{1,p_1}}),\ldots, \mu^{\Lie}_{p_s}(X_{i_{s,1}},\ldots ,X_{i_{s,p_s}}),X_1\bo]_{\mathrm L},\notag
\end{align}
where the generating function of the coefficients $\beta_s$ is
\[\sum_{s=0}^\infty\beta_sx^s= \beta(x)=\frac{x}{\mathrm e^x-1};\]
\begin{align}
&\mu^{\Lie}_n(X_1,\ldots,X_n)=\tag{R}\plabel{R}\\
&= \sum_{\substack{J_1\dot\cup\ldots \dot\cup J_r=\{1,\ldots,n-1\}\\J_l=\{j_{l,1},\ldots,j_{l,q_l}\}\neq
\emptyset \\j_{l,1}<\ldots<j_{l,q_l}}} \tilde \beta_r\cdot \bo[ \mu^{\Lie}_{q_1}(X_{j_{1,1}},\ldots ,X_{j_{1,q_1}}),\ldots\mu^{\Lie}_{q_r}(X_{j_{r,1}},\ldots ,X_{j_{r,q_r}}),X_n\bo]_{\mathrm L},\notag
\end{align}
where the generating function of the coefficients $\tilde \beta_r$ is
\[\sum_{r=0}^\infty\tilde \beta_rx^r= \beta(-x)=\frac{-x}{\mathrm e^{-x}-1};\]
\begin{align}\mu^{\Lie}_1(X_1)=\,X_1\text{ and}&\tag{C}\plabel{C}\\
\mu^{\Lie}_n(X_1,\ldots,X_n)=&\, \sum_{\substack{I_1\dot\cup\ldots \dot\cup I_s\dot\cup  J_1\dot\cup\ldots \dot\cup J_r=\{2,\ldots,n-1\}\\I_k=\{i_{k,1},\ldots,i_{k,p_k}\}\neq
\emptyset, J_l=\{j_{l,1},\ldots,j_{l,q_l}\}\neq
\emptyset \\
i_{k,1}<\ldots<i_{k,p_k}
,j_{l,1}<\ldots<j_{l,q_l}}} \alpha_{s,r}\cdot\notag\\
&\qquad\bo[\quad\bo[\mu^{\Lie}_{p_1}(X_{i_{1,1}},\ldots ,X_{i_{1,p_1}}),\ldots, \mu^{\Lie}_{p_s}(X_{i_{s,1}},\ldots ,X_{i_{s,p_s}}),X_1\bo]_{\mathrm L},\notag\\
&\qquad \quad\,\,\bo[ \mu^{\Lie}_{q_1}(X_{j_{1,1}},\ldots ,X_{j_{1,q_1}}),\ldots,\mu^{\Lie}_{q_r}(X_{j_{r,1}},\ldots ,X_{j_{r,q_r}}),X_n\bo]_{\mathrm L} \quad \bo]\notag
\end{align}
for $n\geq 2$, where the generating function of the coefficients $\alpha_{s,r}$ is
\begin{align}
\sum_{s=0}^\infty\sum_{r=0}^\infty \alpha_{s,r} x^sy^r= \alpha(x,y)&=\frac{\beta(-x-y)-\beta(-y)}{x}\beta(x)\tag{g}\label{eq:cg}\\
&=-\frac{\beta(x+y)-\beta(x)}{y}\beta(-y)\notag\\
&=
{\frac {  -x{{\rm e}^{y}}-{{\rm e}^{-
x}}y+x+y  }{ \left( {{\rm e}^{-y}}-{{\rm e}^{x}} \right)
 \left({{\rm e}^{y}} -1\right)  \left({{\rm e}^{-x}}-1 \right) }}.\notag
\end{align}
\begin{proof}[Note]
Here $\alpha(x,y)$ exists primarily as the RHS of \eqref{eq:cg} line1 or line2, where it is easy to see that
 one has formal power series.
Then expansion to the RHS of \eqref{eq:cg} line3 ($\mathbb Q[[x,y]]$ being zero-divisor free)
 proves equality between the two definitions.
\qedno
\end{proof}
\begin{proof} It is easy to see that $\beta(x)$  has only rational coefficients.
We have here three different recursive definitions for $\mu^{\Lie}_n$, we have to show that they  give the same Lie-polynomials.

The three definitions agree for $n=1$.
By induction, assume that the $\mu^{\Lie}_m$ are well-defined for $m<n$, $n\geq2$.
Consider first the definition of $\mu^{\Lie}_n(X_1,\ldots,X_n)$ according to \eqref{R}. 
This gives
\[
\mu^{\Lie}_n(X_1,\ldots,X_n)_{(R)}= \sum_r \sum_{(X_1,\ldots,X_{n-1} \text{ loc. incr.})}\tilde \beta_r\cdot
\bo[ \underbrace{\mu(\lldots),\ldots,\mu(\lldots)}_{\text{$r$ many $\mu$}},X_n\bo]_{\mathrm L},\notag
\]
where we do not fill out the variables and indices in the $\mu(\lldots)$,
but just note that we have to sum for all locally increasing deployments of the variables $X_1,\ldots,X_{n-1}$.

In each summand, we consider the $\mu(\lldots)$ containing $X_1$, and expand it using \eqref{L}.
(We can do this according to the induction hypothesis).
It yields
\[
= \sum_{r,p,s}  \sum_{(X_2,\ldots,X_{n-1} \text{ loc. incr.})}\tilde \beta_r\beta_s\cdot
\bo[ \underbrace{\mu(\lldots),\ldots,\mu(\lldots)}_{\text{$p$ many $\mu$}},
 \bo[ \underbrace{\mu(\lldots),\ldots,\mu(\lldots)}_{\text{$s$ many $\mu$}},X_1\bo]_{\mathrm L},
\underbrace{\mu(\lldots),\ldots,\mu(\lldots)}_{\text{$r-1-p$ many $\mu$}}
,X_n\bo]_{\mathrm L}.\notag
\]
There we have  a commutator $\bo[\bo[\underbrace{\mu(\lldots),\ldots,\mu(\lldots)}_{\text{$s$ many $\mu$}},X_1\bo]_{\mathrm L},
\bo[\underbrace{\mu(\lldots),\ldots,\mu(\lldots)}_{\text{$r-1-p$ many $\mu$}}
,X_n\bo]_{\mathrm L}\bo]$, which is commutated further by $p$ many $\mu(\lldots)$'s.
Let us distribute those, using the Leibniz rule, between the two terms of the commutators.
We obtain
\[= \sum_{\bar s,\bar r}  \sum_{(X_2,\ldots,X_{n-1} \text{ loc. incr.})}\tilde \alpha_{\bar s,\bar r}\cdot \bo[
 \bo[ \underbrace{\mu(\lldots),\ldots,\mu(\lldots)}_{\text{$\bar s$ many $\mu$}},X_1\bo]_{\mathrm L},\bo[
\underbrace{\mu(\lldots),\ldots,\mu(\lldots)}_{\text{$\bar r$ many $\mu$}}
,X_n\bo]_{\mathrm L}\bo].\notag\]
This is a formally deterministic process, which gives nonzero contributions only for $\bar s+\bar r\leq n-2$
 (because there are only $n-2$ many variables to distribute).
According to our specific method, if $k\leq n-2$, then
\begin{multline}\notag
\sum_{\bar s+\bar r=k}\tilde\alpha_{\bar s,\bar r}x^{\bar s}y^{\bar r}
=\sum_{\substack{r-1,s\geq 0,\\ (r-1)+s=k}}\tilde \beta_r\left(y^{r-1}+\ldots+(x+y)^p y^{r-1-p}+\ldots (x+y)^{r-1}\right) \beta_s x^s\\
 =\sum_{\substack{r-1,s\geq 0,\\ (r-1)+s=k}} \frac{\tilde \beta_r (x+y)^r- \tilde \beta_r y^r}{x}\beta_s x^s=\left(\frac{\beta(-x-y)-\beta(-y)}{x}\beta(x)\right)_{\substack{\text{$k$-homogeneous}\\\text{part in $x,y$}}}.
\end{multline}
We see that our manipulations yield $\tilde\alpha_{\bar s,\bar r}=\alpha_{\bar s,\bar r}$ for $\bar s+\bar r\leq n-2$,
which implies that the definitions of $\mu^{\Lie}_n(X_1,\ldots,X_n)$ according to \eqref{R} and \eqref{C} are the same
(again, terms with $\bar s+\bar r>n-2$ do not appear in either side, as there are only $n-2$ variables to distribute).

The argument that \eqref{L} and \eqref{C} give the same polynomials is analogous.
\end{proof}
\end{propdef}
\snewpage
\begin{commenty}
\begin{remark}
Some values of $\beta_s$ are given by
\[  \begin {array}{c|cccccc}
\beta_{s}&s=0&1&2&3&4&5\\ \hline\\
&1&-\dfrac12&\dfrac {1}{12}&0&-\dfrac {1}{720}&0\,.\\ \end {array}
\]
($B_s=s!\beta_s$ are the Bernoulli numbers.)

Some values of $\alpha_{s,r}$ are given by
\[  \begin {array}{c|ccccc}
\alpha_{s,r}&r=0&1&2&3&4\\ \hline\\
s=0&\dfrac12&\dfrac16&0&-{\dfrac {1}{180}}&0\\ \\
1&-\dfrac16&-\dfrac1{12}&-{\dfrac {1}{120}}&{\dfrac {1}{360}}&{\dfrac {1}{2016}}\\ \\
2&0&{\dfrac {1}{120}}&{\dfrac {1}{240}}&{\dfrac {1}{5040}}&-{\dfrac {1}{4032}}\\ \\
3&{\dfrac {1}{180}}&{\dfrac {1}{360}}&-{\dfrac {1}{5040}}&-{\dfrac {1}{3024}}&-{\dfrac{1}{60480}}\\ \\
4&0&-{\dfrac {1}{2016}}&-{\dfrac {1}{4032}}&{\dfrac {1}{60480}}&{\dfrac {1}{34560}}\,.\end {array}
\]
(The numerators are more complicated for higher indices in both cases.)\qedremark
\end{remark}
\end{commenty}

\begin{prop}\label{p:van} The Lie-polynomials $\mu^{\Lie}_n(X_1,\ldots,X_n)$, $n\geq1$, satisfy $(\mu1)$, $(\mu2)$, $(\mu3)$.

\begin{proof}
Only $(\mu3)$ is nontrivial; so let us consider \eqref{eq:muid}.
It is easy to check that $\mu^{\Lie}_2(X_1,X_2)=\frac12\bo[X_1,X_2\bo]$, which shows the statement for $n=2$.
If $\ldots_2$ is non-empty, then expand $\mu^{\Lie}_n(\ldots_1,X_{k-1},X_k,\ldots_2)-\mu^{\Lie}_n(\ldots_1,X_k,X_{k-1},\ldots_2)$ 
 according to the \eqref{R}-expansions of the $\mu^{\Lie}_n$.
Most of the terms cancel each other except those which contain $X_{k-1}$ and $X_k$ immediately next to each other.
But then the induction hypothesis can be applied to show that it yields the  
 \eqref{R}-expansion of $\mu^{\Lie}_{n-1}(\ldots_1,\bo[X_{k-1},X_k\bo],\ldots_2)$.
If $\ldots_1$ is non-empty, then the \eqref{L}-expansion can be used to prove the identity in the same manner.
\end{proof}
\end{prop}

\snewpage
\section{The existence of $\mu^{\Lie}$ II}\plabel{sec:mu2}

We define a  Lie-permutation $Ii$ of $\{1,\ldots,n\}$ as the following data.
It is a partition   $I_1\dot\cup\ldots \dot\cup I_s=\{1,\ldots,n\}$ such that $\max I_1<\ldots<\max I_s$,
 and a collection of finite sequences $i_{k,1},\ldots,i_{k,p_k}$ 
 such that $\{i_{k,1},\ldots,i_{k,p_k}\}=I_k$, $p_k=|I_k|$ and $i_{k,p_k}=\max I_k$.
\begin{lemma}\plabel{lem:deco}
 The number of Lie-permutations of $\{1,\ldots,n\}$, $n\geq0$, is $n!$.
\begin{proof}
For any Lie-permutation $Ii$ write down the sequence
\[\underbrace{i_{s,1},\ldots,i_{s,p_s}}_{\text{from }I_s},\ldots,\underbrace{i_{k,1},\ldots,i_{k,p_k}}_{\text{from }I_k},\ldots,\underbrace{i_{1,1},\ldots,i_{1,p_1}}_{\text{from }I_1}.\]
This yields a permutation of $\{1,\ldots,n\}$.
From this permutation the Lie-permutation can be reconstructed.
Indeed, in the permutation sequence, the first couple of elements up to `$n$' form
the last partition set $I_s$ with ordering. Then, from the rest, the first couple of elements up to the
maximal element form the   partition set $I_{s-1}$ with ordering; etc.
It is easy to see that we have a bijection between permutations and Lie-permutations.
\end{proof}
\end{lemma}
In what follows let $\mathbb QX_{\Sigma_n}$ be the vector space spanned by the noncommutative monomials
$X_{\sigma(1)}\ldots X_{\sigma(n)}$ in the corresponding noncommutative polynomial ring over $\mathbb Q$.
\begin{prop}\plabel{prop:deco}
Any element of $\mathbb QX_{\Sigma_n}$ can uniquely be written in the form
\begin{equation}
\sum_{Ii\text{ is a Lie-permutation of }\{1,\ldots,n\}} a_{Ii}
[X_{i_{1,1}},\ldots,X_{i_{1,p_1}}]_{\mathrm L}\cdot_\Sigma\ldots \cdot_\Sigma [X_{i_{s,1}},\ldots,X_{i_{s,p_s}}]_{\mathrm L}
\plabel{eq:decoe}
\end{equation}
where $a_{Ii}\in\mathbb Q$. (Here we used ordinary commutators and symmetrized products.)
\begin{proof}
Existence is a consequence of the standard symmetrization argument but applied in the non-commutative polynomial algebra.
This proves that any element is a sum symmetric products of commutator monomials.
Commutator monomials, on the other hand, can be brought into standard form (highest indices on the right in left-iterated commutators).
Uniqueness follows from dimensional reasons, as the number of Lie-partitions of $\{1,\ldots,n\}$ is $n!$, the same
as the dimension of $\mathbb QX_{\Sigma_n}$.
\end{proof}
\end{prop}
\begin{cor}\plabel{cor:deco}
$\mathbb QX_{\Sigma_n}$ allows a direct sum decomposition according to the number of components in the symmetric products in
\eqref{eq:decoe}.

This decomposition is invariant to relabelings
$\ldots,X_k,X_{k+1},\ldots\rightsquigarrow\ldots,X_{k+1},X_k,\ldots$
or substitutions $X_1,\ldots,X_k,\ldots,X_{n-1}\rightsquigarrow X_1,\ldots,[X_k,X_{k+1}],\ldots,X_{n}$.
\begin{proof}
The decomposition property follows from the unicity of \eqref{eq:decoe}.
If we make the substitutions indicated, then the $s$-ply symmetric products go $s$-ply symmetric products,
and they remain so even if the substituted commutator monomials are expanded to standardized ones
(highest indices on the right in left-iterated commutators).
\end{proof}
\end{cor}
(Remark:  A much more effective form of the corollary above is the eigendecomposition along
 the Friedrichs co-shuffle map, cf. \cite{LL1}, and the discussion later.)

Let us  write the monomial $X_1\ldots X_n$ into a form like above:
\begin{equation}
X_1\ldots X_n=\!\!\sum_{Ii\text{ is a Lie-permutation of }\{1,\ldots,n\}}\!\! b_{Ii}
[X_{i_{1,1}},\ldots,X_{i_{1,p_1}}]_{\mathrm L}\cdot_\Sigma\ldots \cdot_\Sigma [X_{i_{s,1}},\ldots,X_{i_{s,p_s}}]_{\mathrm L};
\plabel{eq:Lieid}
\end{equation}
the $b_{Ii}$ are concrete rational numbers.
\begin{defi}
Then let us define
\[\mu^{\Lie}_n(X_1,\ldots, X_n)=\sum_{Ii\text{ is a Lie-permutation of }\{1,\ldots,n\},\text{ of one single block}} b_{Ii}
 \bo[X_{i_{1,1}},\ldots,X_{i_{1,p_1}}\bo]_{\mathrm L},
\]
where we now use Lie-commutators instead of commutators.
\end{defi}

\begin{prop}The Lie-polynomials $\mu^{\Lie}_n(X_1,\ldots,X_n)$, $n\geq1$, satisfy $(\mu1)$, $(\mu2)$, $(\mu3)$.
\begin{proof} Only $(\mu3)$ is nontrivial.
Using Lie algebra rules, both sides of \eqref{eq:muid} can be brought into form
\[\sum_{\chi\in\Sigma_{n-1} }c_{\chi}^{\mathrm{LHS}}\bo[X_{\chi_1},\ldots,X_{\chi_{n-1}},X_n \bo]_{\mathrm L}
\qquad\text{and}\qquad
\sum_{\chi\in\Sigma_{n-1} }c_{\chi}^{\mathrm{RHS}}\bo[X_{\chi_1},\ldots,X_{\chi_{n-1}},X_n \bo]_{\mathrm L}\]
respectively. As the commutator evaluations  would  distinguish them by the $X_n$-ending associative monomials,
it is sufficient to prove equality in the commutator expansion.
Let us consider
\begin{equation}
\underbrace{\ldots_1X_{k-1}X_k\ldots_2}_{n\text{ terms}}
-\underbrace{\ldots_1X_kX_{k-1}\ldots_2}_{n\text{ terms}}
=\underbrace{\ldots_1[X_{k-1},X_k]\ldots_2}_{n-1\text{ terms}}.
\plabel{eq:daor}
\end{equation}
The component of $\ldots_1X_{k-1}X_k\ldots_2$ of symmetric degree $1$ is
$\mu^{\ass}_n(\ldots_1,X_{k-1},X_k,\ldots_2)$ (the commutator evaluation of $\mu^{\Lie}_n(\ldots_1,X_{k-1},X_k,\ldots_2)$).
By Corollary \ref{cor:deco}, the components of symmetric degree $1$ for the other terms in equation \eqref{eq:daor}
are $\mu^{\ass}_n(\ldots_1,X_k,X_{k-1},\ldots_2)$ and $\mu^{\ass}_n(\ldots_1,[X_{k-1},X_k],\ldots_2)$, repectively.
This implies
\[\mu^{\ass}_n(\ldots_1,X_{k-1},X_k,\ldots_2)-\mu^{\ass}_n(\ldots_1,X_k,X_{k-1},\ldots_2)=\mu^{\ass}_n(\ldots_1,[X_{k-1},X_k],\ldots_2),\]
what we wanted to show.
\end{proof}
\end{prop}

\snewpage
\section{From $\mu^{\Lie}$ to the symmetric PBW theorem}\plabel{sec:symm}

\begin{defi}
We define the map $\mu_\Sigma:\bigotimes\mathfrak g\rightarrow  \bigotimes_\Sigma \mathfrak g$
such that
\[
\mu_\Sigma(x_1\otimes\ldots\otimes x_n)=\sum_{\substack{I_1\dot\cup\ldots \dot\cup I_s=\{1,\ldots,n\}\\I_k=\{i_{k,1},\ldots,i_{k,p_k}\}\neq
\emptyset \\i_{k,1}<\ldots<i_{k,p_k}}}\frac1{s!}\,
 \mu^{\Lie}_{p_1}(x_{i_{1,1}},\ldots ,x_{i_{1,p_1}})\otimes\ldots\otimes \mu^{\Lie}_{p_s}(x_{i_{s,1}},\ldots ,x_{i_{s,p_s}})\notag
\]
(and it acts trivially in the $0$th order).
\end{defi}
\begin{propdef}  $\mu_\Sigma:\bigotimes\mathfrak g\rightarrow  \bigotimes_\Sigma \mathfrak g$
descends to a map $\bo\mu_\Sigma:\mathcal U\mathfrak g\rightarrow  \bigotimes_\Sigma \mathfrak g$.
\begin{proof}
It is sufficient to check that
\[\mu_\Sigma\left(\ldots_1\otimes X_{k-1}\otimes X_k\otimes \ldots_2
-\ldots_1\otimes X_{k-1}\otimes X_k\otimes \ldots_2
 -\ldots_1\otimes \bo[X_{k-1}, X_k\bo]\otimes \ldots_2\right)=0,\]
i.~e.~$\mu_\Sigma$ vanishes on the ideal generated by the elements $X\otimes Y-Y\otimes X-\bo[X,Y\bo] $.
This vanishing, when expanded, however, is a consequence of identities \eqref{eq:muid}.
\end{proof}
\end{propdef}
\begin{lemma}\plabel{p:per} Suppose that  $P(X_1,\ldots,X_n)$, with $n\geq 2$,
is a combination of Lie-monomials,
such that in every Lie-monomial every variable appears exactly once. Then
\[\sum_{\sigma\in\Sigma_n}P(X_{\sigma(1)},\ldots,X_{\sigma(n)})=0. \]

In particular, it holds that
\begin{equation}
\sum_{\sigma\in\Sigma_n} \mu^{\Lie}_n(X_{\sigma(1)},\ldots,X_{\sigma(n)})=0 \qquad\text{for}\qquad n\geq2.
\tag{$\mu2^{\mathrm{weak}}$}
\plabel{eq:musymm}
\end{equation}
\begin{proof} This is sufficient to prove for Lie-monomials of $X_1,\ldots,X_n$.
If $P$ is a non-trivial monomial, then it contains an inner Lie-commutator $\bo[X_k,X_l\bo]$, $k\neq l$.
Now, the permutations from $\Sigma_n$ come in pairs $\sigma\in A_n$ and $\sigma\circ(k\,\,l)\in\Sigma_n\setminus A_n$, which
cancel each other in the permuted monomial.
\end{proof}
\end{lemma}
\begin{prop}
 $\bo\mu_\Sigma$ inverts $\bo m_\Sigma$.
\begin{proof} Then
\begin{multline}
\bo\mu_\Sigma\left(\bo m_\Sigma\left(\frac1{n!}\sum_{\sigma\in\Sigma_n}g_{\sigma(1)}\otimes\ldots\otimes g_{\sigma(n)}\right)\right)=
\mu_\Sigma\left(\frac1{n!}\sum_{\sigma\in\Sigma_n}g_{\sigma(1)}\otimes\ldots\otimes g_{\sigma(n)}\right)\\
=\frac1{n!}\sum_{\sigma\in\Sigma_n}\mu^{\Lie}_1(g_{\sigma(1)})\otimes_\Sigma\ldots\otimes_\Sigma \mu^{\Lie}_1(g_{\sigma(n)})=\frac1{n!}\sum_{\sigma\in\Sigma_n}g_{\sigma(1)}\otimes\ldots\otimes g_{\sigma(n)}.
\notag
\end{multline}
Indeed,  according to the previous definition,
 $\mu_\Sigma=\bo\mu_\Sigma\circ\bo m$, furthermore $\bo m_\Sigma$ is just a restriction of $\bo m$;
 this implies the first equality.
The second equality is due to the fact that the higher $\mu_h$ ($h\geq2$) vanish under symmetrization
(Lemma \ref{p:per}).
The third one is true due to $\mu^{\Lie}_1(X_1)=X_1$ and that
the symmetrization of the symmetrization is the symmetrization. This proves
$\bo\mu_\Sigma\circ\bo m_\Sigma=\id_{\bigotimes_\Sigma\mathfrak g}$.
In particular, $\bo m_\Sigma$ is injective.
Then, the surjectivity of $\bo m_\Sigma$ implies bijectivity, and, in fact, the inverse relationship.
\end{proof}
\end{prop}

This, in particular, proves the symmetric global version of the PBW theorem, i.~e.~that $\bo m_\Sigma$ is an isomorphism.
We have not used any unicity result for $\mu^{\Lie}_n$ but just existence.
\snewpage

\textbf{Note on unicity.} One can see from the proof that the requirement  that the $\mu^{\Lie}_n$
are Lie-polynomials satisfying $(\mu1)$,  $(\mu2)$, $(\mu3)$
can be replaced by the weaker requirement that the $\mu^{\Lie}_n$
are multilinear functions $\mathfrak g^n\rightarrow\mathfrak g$ satisfying $(\mu1^*)$,  $(\mu2^*)$, $(\mu3^*)$,
where $(\mu1^*)$ corresponds to $(\mu1)$,  $(\mu2^*)$ corresponds to \eqref{eq:musymm},
$(\mu3^*)$ corresponds to $(\mu3)$.

As $\boldsymbol\mu_\Sigma$ inverts $\boldsymbol m_\Sigma$, and the inverse is unique, in particular
its part of order $1$ is unique, this demonstrates the uniqueness of the $\mu^{\Lie}_n$ as Lie algebraic functions.
But then this also shows that their effect are the same as given by the Lie-polynomials we have constructed.
Applied to free Lie algebras, it also follows that the constructions of $\mu^{\Lie}_n$ in the previous two
sections are the same.
However,

\snewpage
\section{More on $\mu^{\Lie}$}\plabel{sec:furthermu1}

The simplest way to show unicity for the $\mu^{\Lie}_n$ as Lie-polynomials
is just to apply a fixed symmetric rearrangements process in the variables.
Due to the $(\mu3)$ and  \eqref{eq:musymm}, and $(\mu1)$, at the end we obtain
$\mu^{\Lie}_n(X_1,\ldots,X_n)=\mu^{\Lie}_1(P(X_1,\ldots,X_n))=P(X_1,\ldots,X_n)$
where $P$ is a Lie-polynomial depending on the rearrangements process but not the $\mu^{\Lie}_n$.

If $\mu^{\Lie}_n(X_1,\ldots,X_n)$ is commutator-evaluated to $\mu^{\ass}_n(X_1,\ldots,X_n)$,
then one finds %, by relatively simple combinatorial arguments,
 that
\begin{align}
\mu_n^{\ass}(X_1,\ldots,X_n) &=\frac{\partial^n}{\partial t_1\ldots \partial t_n}\log(\exp(t_1 X_1)\ldots \exp(t_n X_n))\biggr|_{t_1,\ldots,t_n=0} \plabel{eq:sE}
\\&=\sum_{\sigma\in\Sigma_n}
\underbrace{(-1)^{\des(\sigma)}\frac{\des(\sigma)!\asc(\sigma)!}{n!}}_{\mu_\sigma:=}  X_{\sigma(1)}\ldots X_{\sigma(n)};
\plabel{eq:sEE}
\end{align}
where $\asc(\sigma)$ denotes the number of ascents,
i.~e.~the number of pairs such that $\sigma(i)<\sigma(i+1)$;
and  $\des(\sigma)$ denotes the number of  descents,
i.~e.~the number of pairs such that $\sigma(i)>\sigma(i+1)$.
This is the result of Dynkin \cite{Dyy} (cf.~Wall \cite{Wa}, Achilles, Bonfiglioli \cite{AB},
 and also Strichartz \cite{St}), Solomon \cite{SS}, or
Mielnik, Pleba\'nski \cite{MP} (with various emphases).

One possible approach for a proof is as follows:
One can see that the non-commutative polynomial expressions $\tilde\mu_n^{\ass}$ involved in 
the three relevant places of \eqref{eq:sE}--\eqref{eq:sEE} satisfy
\begin{enumerate}
\item[($\mu1'$)] $\tilde\mu^{\ass}_1(X_1)=X_1$;
\item[($\mu2'$)] $\tilde\mu_n^{\ass}(X,\ldots,X)=0$ holds for $n\geq2$, but $\tilde\mu_n^{\ass}(X_{1},\ldots,X_{n})$ is linear in its variables;
\item[($\mu3'$)] the identities  \eqref{eq:muid} hold but meant with $\tilde\mu_n^{\ass}$ and ordinary commutators.
\end{enumerate}
(Very moderate combinatorial arguments are sufficient here with respect to RHS\eqref{eq:sE} and RHS\eqref{eq:sEE}, cf. \cite{LL1}.)
Then one can argue that   ($\mu1'$)--($\mu3'$) determine the $\tilde\mu_n^{\ass}$ in the associative
setting, again by using the symmetric rearrangement procedure.

\textbf{From $\tilde\mu_n^{\ass}(X,\ldots,X)$ to $\mu^{\Lie}$.}
At first sight it might be disheartening that the notable expressions RHS\eqref{eq:sE}
and RHS\eqref{eq:sEE} did not lead directly to the existence of Lie algebraic $\mu_n^{\Lie}$.
But they do, even if in not an explicit way:
 one can show that the existence of $\tilde\mu_n^{\ass}$ implies the existence of $\mu_n^{\Lie}$.
Indeed, let us just try the symmetric rearrangement process.
We cannot be sure about its Lie algebraic
consistency, but we can record our experiences in terms of free non-associative algebras (magmas, brace algebras).
In that way, we obtain
\[X_1\otimes\ldots\otimes X_n\rightsquigarrow H^{(n)}_n+\ldots+H^{(1)}_n,\]
where $H^{(i)}_n$ is an $i$-ply symmetric tensor over the nonassociative algebra $ \mathrm F_{\mathbb Q}^{\nonass}[X_1,\ldots,X_n]$
with overall homogeneity degree $1,\ldots,1$ in its variables $X_1,\ldots,X_n$.
We spell out that $H^{(1)}_n=P^{\nonass}_n(X_1,\ldots,X_n)\in\mathrm F_{\mathbb Q}^{\nonass}[X_1,\ldots,X_n]$.
We can resolve $P^{\nonass}_n(X_1,\ldots,X_n)$ to the Lie-polynomial $P^{\Lie}_n(X_1,\ldots,X_n)$, which, by commutator-evaluation, can be resolved further to
$P^{\ass}_n(X_1,\ldots,X_n)$.
By the associative unicity we know that $P^{\ass}_n(X_1,\ldots,X_n)=\tilde\mu_n^{\ass}(X_1,\ldots,X_n)$.
As uniformy $1$-homogeneous free Lie algebra elements can be brought to shape
$\sum_{\chi\in\Sigma_{n-1} }c_{\chi} \bo[X_{\chi_1},\ldots,X_{\chi_{n-1}},X_n \bo]_{\mathrm L}$,
where the commutator-evaluation distinguishes (by the $X_n$-ending monomials),
the commutator evaluation is faithful here,
thus
there no is choice for $P^{\Lie}_n$ but to satisfy $(\mu3)$  with $(\mu1)$ and $(\mu2)$ being trivial.
Therefore the $P^{\Lie}_n$ give $\mu_n^{\Lie}$.
In that way, we obtain ``The existence of $\mu^{\Lie}$ III'' with respect to
RHS\eqref{eq:sE} and $\mu^{\Lie}$  IV'' with respect to
RHS\eqref{eq:sEE}.

\textbf{Explicit forms (from $\mu^{\ass}$-IV).}
In the following arguments the actual shape of \eqref{eq:sEE} plays no role, thus we just write
\begin{equation}
\mu^{\ass}_n(X_1,\ldots,X_n)=\sum_{\sigma\in\Sigma_n}\mu_\sigma X_{\sigma(1)}\ldots X_{\sigma(n)}.
\plabel{eq:EE}
\end{equation}

Let us fix an arbitrary element $k\in\{1,\ldots,n\}$.
Now, $\mu^{\Lie}_n(X_1,\ldots,X_n)$ is a Lie-polynomial, thus, using standard commutator rules,
we can write it as linear combination of terms $\bo[X_{i_1},\ldots,X_{i_{n-1}} ,X_k\bo]_{\mathrm L}$,
where $\{i_1,\ldots,i_{n-1}\}=\{1,\ldots,n\}\setminus\{k\}$. However, evaluated in the
 noncommutative polynomial algebra, such a commutator expression gives
only one monomial contribution $X_{i_1}\ldots X_{i_{n-1}} X_k $ such that the last term is $X_k$.
Thus, the coefficient of $\bo[X_{i_1},\ldots,X_{i_{n-1}} ,X_k\bo]_{\mathrm L}$ can be read off from \eqref{eq:EE}.
We find that
\begin{equation}
\mu_n^{\Lie}(X_1,\ldots,X_n)=
\sum_{\sigma\in\Sigma_n, \sigma(n)=k}\mu_\sigma \bo[X_{\sigma(1)},\ldots ,X_{\sigma(n-1)},X_k\bo]_{\mathrm L}.\plabel{eq:EEk}
\end{equation}
Cf.~ Dynkin \cite{Dyy} and Arnal, Casas, Chiralt \cite{ACC}.
Averaging \eqref{eq:EEk} for all possible $k$, we obtain
\begin{equation}
\mu_n^{\Lie}(X_1,\ldots,X_n)=\frac1n
\sum_{\sigma\in\Sigma_n}\mu_\sigma \bo[X_{\sigma(1)},\ldots ,X_{\sigma(n-1)},X_{\sigma(n)}\bo]_{\mathrm L},\plabel{eq:EEn}
\end{equation}
which is the standard `Dynkinized' form of \eqref{eq:EE}, cf.~the above mentioned sources.

%\begin{commenty}
In the arguments above we have not used the PBW theorem or any nontrivial fact about free Lie algebras,
 but only the fact that the commutator-evaluation is sufficiently restrictive here.
Although these arguments are of Dynkin \cite{Dyy},
they subsume the standard form of the Dynkin--Specht--Wever lemma.
(Nevertheless, \eqref{eq:EEn} can be thought as an application of the  Dynkin--Specht--Wever lemma,
and \eqref{eq:EEk} can be thought as an application of a weighted  Dynkin--Specht--Wever lemma.)
%\end{commenty}

\textbf{Connection to the BCH expansion (from $\mu^{\ass}$-III).} This is most well-known,
\begin{align*}
\BCH^{\ass}_{n}(X,Y)&\equiv\sum_{r=0}^n\log\left((\exp X)(\exp Y)\right)_{\substack{
\text{$r$-homogeneous part in $X$} \\ \text{$n-r$-homogeneous part in $Y$}}}\\
&=\sum_{r=0}^n\frac1{r!}\frac1{(n-r)!} \frac{\partial^r}{\partial t^r}
\frac{\partial^{n-r}}{\partial \tau^{n-r}}\log(\exp(t X)  \exp( \tau Y))\biggr|_{t,\tau=0}\\
&=\sum_{r=0}^n
\frac1{r!}\frac1{(n-r)!} \frac{\partial^r}{\partial t_1\ldots \partial t_r}
\frac{\partial^{n-r}}{\partial \tau_1\ldots \partial \tau_{n-r}}\\&\qquad\log(\exp((t_1+\ldots+t_r) X)  \exp( (\tau_1+\ldots+\tau_{n-r})Y))\biggr|_{t_1,\ldots,\tau_{n-r}=0}\\
&=\sum_{r=0}^n\frac1{r!}\frac1{(n-r)!}
\mu^{\ass}_{n}(\underbrace{X,\ldots,X}_{r\text{ times}},\underbrace{Y,\ldots, Y}_{n-r\text{ times}}).
\end{align*}
Thus, a Lie algebraic lift is provided by
\begin{equation}
\BCH^{\Lie}_n(X,Y)=\sum_{r=0}^n\frac1{r!}\frac1{(n-r)!}
\mu^{\Lie}_{n}(\underbrace{X,\ldots,X}_{r\text{ times}},\underbrace{Y,\ldots, Y}_{n-r\text{ times}}).
\plabel{eq:mubchlie}
\end{equation}
The existence of such a lift also follows from Fridriechs criterion (which, in turn, can
also be thought as yet another appearance of the symmetric rearrangement method, cf. \cite{LL1}.)
The uniqueness of the lift, in general, follows from the full Dynkin--Specht--Wever lemma
 or from the PBW theorem (the point is that the commutator representation of free Lie algebras is faithful,
 originally a theorem of Magnus \cite{MM} and Witt \cite{W}).
The Magnus commutator can also be constructed from the BCH expression but for our purposes something like
 ($\mu1'$)--($\mu3'$) should be proven anyway, thus this latter approach would be less economical for us.

\begin{remark}
While \eqref{eq:EEk} is sufficiently economical (i.~e.~redundancy-free) with respect to the Magnus expansion,
it is not so in terms of the BCH specialization \eqref{eq:mubchlie}.
However, an essentially redundancy-free ``explicit'' form was obtained by Macdonald \cite{Mac} using the classical Hall bases
of M. Hall \cite{Ha} (whose work bases on P. Hall \cite{Ha0} and Magnus \cite{MM}).
Macdonald \cite{Mac} builds on the ideas of Dynkin \cite{Dyy} and Goldberg \cite{Gol} in greater depth.
For generalization, see Bandiera, Simi \cite{BS}.
\qedremark
\end{remark}

\textbf{Expressions involving  coproducts (from $\mu^{\ass}$-II).}
The idea of Section \ref{sec:mu2} is decomposition by symmetric components.
Technically this can be realized by co-shuffles.
In the noncommutative polynomial algebra let $\Delta$ be the coproduct, which is the substitution operation sending
 the variable $X_\lambda$ to $X_\lambda\otimes1+1\otimes X_\lambda$, and let $m$ be the multiplication operation
 sending $S_1\otimes S_2$ to $S_1S_2$.
One defines the Friedrichs co-shuffle map as $F_{\langle2\rangle}=m\cdot\Delta$.
It is the observation of Friedrichs \cite{Fri} that for a commutator polynomial $C$ one has
 $\Delta(C)=C\otimes1+1\otimes C$, therefore $F_{\langle2\rangle}(C)=2C$.
One can easily see that if $S$ is a $s$-ply symmetrized product of commutator monomials, then $F_{\langle2\rangle}(S)=2^sS$.
Using $F_{\langle2\rangle}$ one can already decompose to eigenspaces (cf. \cite{LL1}) but here we proceed differently.
More generally, one can define the $p$-ply co-shuffle map as $F_{\langle p\rangle}=m^{p-1}\cdot\Delta^{p-1}$
 (associativity and coassociativity used).
One can similarly see that if $S$ is a $s$-ply symmetrized product of commutator monomials, then $F_{\langle p\rangle}(S)=p^sS$.
($p=1$ corresponds to the identity, $p=0$ corresponds to the counit).
For $p\geq1$, we may consider
\begin{equation}
F^{(1)}_{\langle p\rangle}=
\left|\begin{matrix}
1&1&\cdots&1\\
1&2&\cdots&2^p\\
\vdots&\vdots&&\vdots\\
1&p&\cdots&p^p\\
\end{matrix}\right|^{-1}
\left|\begin{matrix}
F_{\langle1\rangle}&1&\cdots&1\\
F_{\langle2\rangle}&2&\cdots&2^p\\
\vdots&\vdots&&\vdots\\
F_{\langle p\rangle}&p&\cdots&p^p\\
\end{matrix}\right|
=\sum_{i=1}^p\frac{(-1)^{i-1}}{i}\binom{p}{i}F_{\langle i\rangle}.
\plabel{eq:pref1}
\end{equation}
Then, by simple linear algebra, for $1\leq n\leq p$,
\begin{equation}
\mu^{\ass}_n(X_1,\ldots, X_n)=F^{(1)}_{\langle p\rangle}(X_1\ldots X_n);
\plabel{eq:pref2}
\end{equation}
the case $n=p$ is practically sufficient to obtain a presentation for $\mu^{\ass}_n(X_1,\ldots, X_n)$.

However, if one already uses the coproduct,
 it is not  necessary to write out the Dynkin--Magnus commutators, but using an other approach,
 the PBW isomorphism can be constructed directly from the coproduct.
That is a connection to the approach of Milnor, Moore \cite{MiMo}.

\textbf{Comparison of constructions.}
$\mu^{\Lie}$-I using the Magnus recursion formulae stands alone as the construction which uses the less,
 although it is quite artificial.
$\mu^{\Lie}$-II, $\mu^{\Lie}$-III, $\mu^{\Lie}$-IV, all
use the polynomial method of Dynkin
 and the symmetric rearrangement process.
With $\mu^{\Lie}$-II, the technical machinery was kept at a low level,
 but it can be made more explicit using the coproduct.
$\mu^{\Lie}$-III gives the most direct approach to the BCH and Magnus expansions.
$\mu^{\Lie}$-IV is the closest one to Solomon \cite{SS}.

\textbf{In stead of $\mu^{\Lie}$-V.} We finish with an observation.
Assume that we know the PBW theorem for the free Lie algebra $\mathrm F_{\mathbb Q}^{\Lie}[X_1,\ldots]$
somehow, by any other method.
By simple universal algebraic arguments $\mathcal U\mathrm F_{\mathbb Q}^{\Lie}[ X_1,\ldots,]$ is naturally isomorphic to the noncommutative
 polynomial algebra $\mathrm F_{\mathbb Q}[X_1,\ldots]$.
Then symmetric degree $1$ part of $(\boldsymbol m_\Sigma)^{-1}(X_1\cdot\ldots\cdot X_n)$
 can naturally be used for $\mu_n^{\Lie}(X_1,\ldots,X_n)$.
In particular, having the PBW theorem for $\mathrm F_{\mathbb Q}^{\Lie}[X_1,\ldots]$
 implies the symmetric PBW theorem in general.

\snewpage
\section{$\mathcal U\mathfrak g$ as a direct construction and result of Nouaz\'e--Revoy type}\plabel{sec:udir}
Still assume $\mathbb Q\subset K$.
Let us define the maps
$\bch_{n,m}:\bodot^n\mathfrak g\otimes\bodot^m\mathfrak g\rightarrow\mathfrak g$,
for $n+m\geq1$, such that
\[\bch_{n,m}(a_1\odot\ldots\odot a_n,b_1\odot\ldots \odot b_m)=\frac1{n!m!}\sum_{
\sigma\in\Sigma_n, \,\chi\in\Sigma_m}
\mu^{\Lie}_{n+m}\left(a_{\sigma(1)},\ldots,a_{\sigma(n)},b_{\chi(1)},\ldots,b_{\chi(m)}  \right).\]

Considering the natural correspondence between
$a_1\odot\ldots\odot a_n$ and $a_1\otimes_\Sigma \ldots\otimes_\Sigma  a_n
=\frac1{n!}a_{\sigma(1)}\otimes\ldots\otimes a_{\sigma(n)}$,
we obtain
\begin{prop}
$\mathcal U\mathfrak g$ is naturally isomorphic to $\bigodot\mathfrak g$
endowed with a product rule $\cdot_{\mathcal U}$ such that
\begin{multline}
a_1\odot\ldots\odot a_n\cdot_{\mathcal U}b_1\odot\ldots \odot b_m =
\sum_{\substack{I_1\dot\cup\ldots \dot\cup I_s=\{1,\ldots,n+m\}\\
I_k=\{i_{k,1},\ldots,i_{k,p_k},n+j_{k,1},\ldots,n+j_{k,q_k}\}\neq
\emptyset \\i_{k,1}<\ldots<i_{k,p_k}\leq n<n+j_{k,1}<\ldots<n+j_{k,q_k}
\\\min I_1<\ldots<\min I_s}}\\\,\,\,\,
\bch_{p_1,q_1}(a_{i_{1,1}}\odot\ldots\odot a_{i_{1,p_1}},
  b_{j_{1,1}}\odot\ldots\odot b_{j_{1,q_1}})
  \odot
  \ldots\odot
 \bch_{p_1,q_1}(a_{i_{s,1}}\odot\ldots\odot a_{i_{s,p_s}},
  b_{j_{s,1}}\odot\ldots\odot b_{j_{s,q_s}}).
 \plabel{eq:diru}
 \end{multline}
\begin{proof}
Indeed, we  have linear isomorphisms $\bo m_\Sigma$ / $\bo \mu_\Sigma$ between
the two modules.
Regarding the product structure, if we resolve $\odot$ as $\otimes_\Sigma $, take the tensor product,
evaluate by
$\bo\mu_\Sigma$, and resolve $\otimes_\Sigma $ back to $\odot$, then we obtain the product rule as above.
\end{proof}
\end{prop}
\snewpage
Thus, a direct construction $\mathcal U_{\mathrm{ dir}} \mathfrak g$ for $\mathcal U\mathfrak g$, in case  $\mathbb Q\subset K$, would
 simply be $\bigodot\mathfrak g$ endowed with the product rule \eqref{eq:diru}.
(Cf.~ Cartier \cite{Ct}.)
Checking well-definedness directly is not particularly hard, but checking the arithmetics for associativity is not that easy.
Nevertheless, we know that the arithmetics works out, because the proposition above holds
 for the free Lie algebra over the rational numbers.

In particular, it works out in the case of the free $k$-nilpotent Lie algebra, where the identity
 $\bo[X_1,\ldots,X_{k+1}\bo]_{\mathrm L}=0$ holds.
In this case, we can consider the evaluator $\bch_{n,m}$, $n+m\geq k+1$ as identically $0$.
In particular, the associativity works out only using $\bch_{n,m}$, $n+m\leq k$ and the $k$-nilpotency rule.
Here,  $\bch_{n,m}$, $n+m\leq k$ can be defined using only the ring $\mathbb Z[\frac1{k!}]$;
 indeed, in the ``symmetric rearrangement procedure'' leading to $\mu^{\Lie}_{n+m}$ we use symmetrizations  up to $k$ elements only,
 and also in the definitions of $\bch_{n,m}$.
(Or see \eqref{eq:sEE}--\eqref{eq:EE}--\eqref{eq:EEk}, cf.~\eqref{eq:pref1}--\eqref{eq:pref2}.)
Now, the free $k$-nilpotent Lie algebra over $\mathbb Z[\frac1{k!}]$ naturally embeds into
 the free $k$-nilpotent Lie algebra over $\mathbb Q$.
In fact, the free $k$-nilpotent Lie algebra (but not its universal enveloping algebra)
 naturally embeds to the $k$-nilpotent noncommutative polynomial algebra by the commutator representation.
Thus, the associativity computation works out in the free $k$-nilpotent Lie algebra over $\mathbb Z[\frac1{k!}]$.
Therefore, it works out in any $k$-nilpotent Lie algebra with $\frac1{k!}\in K$.
Thus, in that case, $\mathcal U_{\mathrm{dir}} \mathfrak g$ yields an associative algebra.
$\mathcal U_{\mathrm{ dir}} \mathfrak g$ is generated by $\mathfrak g$,
 thus we have a natural factorization map $\mathcal U \mathfrak g\rightarrow  \mathcal U_{\mathrm{ dir}} \mathfrak g$.
Regarding the filtration induced by the image of $\botimes^n\mathfrak g$,
 this induces a natural factor map $ \bodot^{n+1}\mathfrak g/Z_n\rightarrow \bodot^{n+1}\mathfrak g$.
This, however, implies that $Z_n$ is $0$.
In particular, we obtain
\begin{theorem}If $\mathfrak g$ is $k$-nilpotent and $1,\ldots,\frac1k\in K$, then

(o) $\mathcal U_{\mathrm{ dir}} \mathfrak g$ can be defined formally;

(a) $\mathcal U\mathfrak g$ is naturally isomorphic to $\mathcal U_{\mathrm{ dir}} \mathfrak g$; and

(b) the (local) PBW theorem holds for $\mathfrak g$.
\begin{proof}
(a) and (b) are both implied by $Z_n=0$.
\end{proof}
\end{theorem}
This is a generalization of the result of Nouaz\'e, Revoy \cite{NR}.

\snewpage

\end{document}